\documentclass[10pt]{amsart}
\usepackage{amssymb,graphicx,latexsym,xypic,ae}
\newcommand{\sym}{\mathcal{S}}
\newcommand{\N}{{\mathbb N}}

\newcommand{\emptyword}{{\epsilon}}

\DeclareMathOperator{\red}{\mathrm{red}}

\theoremstyle{plain}
\newtheorem{theorem}{Theorem}
\newtheorem*{theorem*}{Theorem}

\newtheorem*{corollary*}{Corollary}
\newtheorem*{lemma*}{Lemma}
\newtheorem{proposition}[theorem]{Proposition}
\newtheorem*{proposition*}{Proposition}

\theoremstyle{definition}
\newtheorem{definition}[theorem]{Definition}
\newtheorem*{definition*}{Definition}

\newtheorem*{example*}{Example}

\theoremstyle{remark}

\newtheorem*{remark*}{Remark}

\newcommand{\bintree}{\xymatrix@!R@!C@=2ex@M=0.3ex}
\newcommand{\ulbintree}{\xymatrix@!R@!C@=1.3ex@M=0ex}
\newcommand{\LEFT}{\ar@{-}[dl]}
\newcommand{\RIGHT}{\ar@{-}[dr]}
\newcommand{\LR}{\LEFT\RIGHT}

\newcommand{\Circ}{\scriptstyle{\bigcirc}}

\def\dd{\makebox[1.2ex]{\rule[.6ex]{.8ex}{.15ex}}}

\def\ab{\ensuremath{12}}
\def\ba{\ensuremath{21}}

\def\axbc{\ensuremath{1{\dd}23}}
\def\axcb{\ensuremath{1{\dd}32}}
\def\bxac{\ensuremath{2{\dd}13}}
\def\bxca{\ensuremath{2{\dd}31}}
\def\cxab{\ensuremath{3{\dd}12}}
\def\cxba{\ensuremath{3{\dd}21}}

\def\abxc{\ensuremath{12{\dd}3}}
\def\acxb{\ensuremath{13{\dd}2}}
\def\baxc{\ensuremath{21{\dd}3}}
\def\bcxa{\ensuremath{23{\dd}1}}
\def\caxb{\ensuremath{31{\dd}2}}
\def\cbxa{\ensuremath{32{\dd}1}}

\def\axcxb{\ensuremath{1{\dd}3{\dd}2}}
\def\bxaxc{\ensuremath{2{\dd}1{\dd}3}}
\def\bxcxa{\ensuremath{2{\dd}3{\dd}1}}
\def\cxaxb{\ensuremath{3{\dd}1{\dd}2}}

\begin{document}

\title[Avoiding a pair of Babson-Steingr\'{\i}msson patterns]{Enumerating permutations avoiding a pair of
  Babson-Steingr\'{\i}msson patterns}

\author{Anders Claesson}
\address{Matematik\\
  Chalmers tekniska h\"ogskola och G\"oteborgs universitet\\
  S-412 96 G\"oteborg, Sweden}
\email{claesson@math.chalmers.se}

\author{Toufik Mansour}
\address{
  Department of Mathematics\\
  Chalmers University of Technology\\
  S-412 96 G\"oteborg, Sweden}
\email{toufik@math.chalmers.se}

\date{\today}

\begin{abstract}
  Babson and Steingr\'{\i}msson introduced
  generalized permutation patterns that allow the requirement that two
  adjacent letters in a pattern must be adjacent in the permutation.
  Subsequently, Claesson presented a complete solution for the number of
  permutations avoiding any single pattern of type
  $(1,2)$ or $(2,1)$. For eight of these twelve patterns the answer is
  given by the Bell numbers. For the remaining four the answer is
  given by the Catalan numbers.

  In the present paper we give a complete solution for the number of
  permutations avoiding a pair of patterns of type $(1,2)$ or $(2,1)$.
  We also conjecture the number of permutations avoiding the patterns
  in any set of three or more such patterns.
\end{abstract}

\maketitle\thispagestyle{empty}

\section{Introduction}

Classically, a pattern is a permutation $\sigma\in\sym_k$, and a
permutation $\pi\in\sym_n$ avoids $\sigma$ if there is no subword of
$\pi$ that is order equivalent to $\sigma$. For example,
$\pi\in\sym_n$ avoids $132$ if there is no $1\leq i < j < k\leq n$
such that $\pi(i) < \pi(k) < \pi(j)$. We denote by $\sym_n(\sigma)$ the set
permutations in $\sym_n$ that avoids $\sigma$.

The earliest result to an instance of finding $|\sym_n(\sigma)|$ seems to be
MacMahon's enumeration of $S_n(123)$, which is implicit in chapter
V of \cite{MacMahon}. The first explicit result seems to be
Hammersley's enumeration of $S_n(321)$ in \cite{Hammersley}. In
\cite[Ch. 2.2.1]{KnuthVol1} and \cite[Ch. 5.1.4]{Kn73v3} Knuth
shows that for any $\sigma \in S_3$, we have $|S_n(\sigma)| =
C_n=\frac{1}{n+1}\binom {2n} n$, the $n$th Catalan number. Later
Simion and Schmidt ~\cite{SiSc85} found the cardinality of
$\sym_n(P)$ for all $P\subseteq \sym_3$.

In \cite{BaSt00} Babson and Steingr\'{\i}msson introduced generalized
permutation patterns that allow the requirement that two adjacent
letters in a pattern must be adjacent in the permutation.  The
motivation for Babson and Steingr\'{\i}msson in introducing these
patterns was the study of Mahonian statistics. Two examples of such
patterns are $\axcb$ and $\acxb$ ($\axcb$ and $\acxb$ are of type
$(1,2)$ and $(2,1)$ respectively). A permutation $\pi = a_1 a_2 \cdots
a_n$ avoids $\axcb$ if there are no subwords $a_i a_j a_{j+1}$ of
$\pi$ such that $a_i<a_{j+1}<a_j$. Similarly $\pi$ avoids $\acxb$ if
there are no subwords $a_i a_{i+1} a_{j}$ of $\pi$ such that
$a_i<a_{j}<a_{i+1}$.

Claesson \cite{Cl01} presented a complete solution for the number of
permutations avoiding any single pattern of type $(1,2)$
or $(2,1)$ as follows.

\begin{proposition}[Claesson \cite{Cl01}]\label{claesson1}
  Let $n\in\N$. We have
  $$
  |\sym_n(p)| =
  \begin{cases}
    B_n & \text{if }\,
    p\in\{
    \axbc, \cxba, \abxc, \cbxa,
    \axcb, \cxab, \baxc, \bcxa
    \},\\
    C_n &\text{if }\,
    p\in\{
    \bxac, \bxca, \acxb, \caxb
    \},\\
  \end{cases}
  $$
  where $B_n$ and $C_n$ are the $n$th Bell
  (\# ways of placing $n$ labelled balls into $n$ indistinguishable
  boxes, see~\cite[A000110]{SP})
  and Catalan numbers, respectively.
\end{proposition}
In addition, Claesson gave some results for the number of permutations
avoiding a pair of patterns.
\begin{proposition}[Claesson \cite{Cl01}]\label{claesson2}
  Let $n\in\N$. We have
  $$
  |\sym_n(\axbc,\,\abxc)| = B^*_n,\,
  |\sym_n(\axbc,\,\axcb)| = I_n,\text{ and }\,
  |\sym_n(\axbc,\,\acxb)| = M_n,
  $$
  where $B^*_n$ is the $n$th Bessel number (\# non-overlapping
  partitions of $[n]$ (see \cite{FlSc90})), $I_n$ is the number of
  involutions in $\sym_n$, and $M_n$ is the $n$th Motzkin number
  (\# ways of drawing any number of nonintersecting chords
  among $n$ points on a circle, see~\cite[A001006]{SP}).
\end{proposition}

This paper is organized as follows. In Section~\ref{prel} we define
the notion of a pattern and some other useful concepts.  For a proof
of Proposition~\ref{claesson1} we could refer the reader to
\cite{Cl01}. We will however prove Proposition~\ref{claesson1} in
Section~\ref{single} in the context of binary trees. The idea being
that this will be a useful aid to understanding of the proofs of
Section~\ref{double}. In Section~\ref{double} we give a solution for
the number of permutations avoiding any given pair of patterns of type
$(1,2)$ or $(2,1)$. These results are summarized in the following
table.
\begin{center}
\quad\quad
$
\begin{array}{|c|c|c|c|c|c|c|c|c|c|c|}
  \hline
  \# \text{ pairs}& 2 & 2& 4& 34& 8& 2&4&4&4&2\\ \hline
  |\sym_n(p,q)| & 0&2(n-1) &\binom n 2 + 1 & 2^{n-1} & M_n & a_n & b_n & I_n & C_n
  & B_n^*\\
  \hline
\end{array}
$
\quad\quad
\mbox{
\begin{minipage}{45ex}
\begin{align*}
  \intertext{Here}
  \sum_{n\geq 0} a_n x^n &=
  \frac 1 {1-x-x^2\sum\limits_{n\geq 0} B^*_n x^n}\\
  \intertext{and}
  b_{n+2} &= b_{n+1} + \sum_{k=0}^n\binom n k b_k.\\
\end{align*}
\end{minipage}
}
\end{center}

Finally, in Section~\ref{multi} we conjecture the sequences
$\{\#\sym_n(P)\}_n$ for sets $P$ of three or more patterns of type $(1,2)$
or $(2,1)$.

\section{Preliminaries}\label{prel}

By an \emph{alphabet} $X$ we mean a non-empty set. An element of $X$
is called a \emph{letter}. A \emph{word} over $X$ is a finite sequence
of letters from $X$. We consider also the \emph{empty word}, that is,
the word with no letters; it is denoted by $\emptyword$. Let
$w=x_1x_2\cdots x_n$ be a word over $X$. We call $|w|:=n$ the
\emph{length} of $w$. A \emph{subword} of $w$ is a word
$v=x_{i_1}x_{i_2}\cdots x_{i_k}$, where $1 \leq i_1<i_2<\cdots<i_k
\leq n$.

Let $[n]:=\{1,2,\ldots,n\}$ (so $[0]=\emptyset$). A \emph{permutation}
of $[n]$ is bijection from $[n]$ to $[n]$. Let $\sym_n$ be the set of
permutations of $[n]$, and $\sym = \cup_{n\geq 0} \sym_n$.  We shall
usually think of a permutation $\pi$ as the word
$\pi(1)\pi(2)\cdots\pi(n)$ over the alphabet $[n]$.

Define the \emph{reverse} of $\pi$ by $\pi^r(i) = \pi(n+1-i)$, and
define the \emph{complement} of $\pi$ by $\pi^c(i) = n+1-\pi(i)$,
where $i\in[n]$.

For each word $w=x_1x_2\cdots x_n$ over the alphabet
$\{1,2,3,4,\ldots\}$ without repeated letters, we define the
\emph{reduction} of $w$, which we denote $\red(w)$, by
$$
\red(w) = a_1 a_2 \cdots a_n \,,\;\text{ where }\; a_i =
|\{j\in [n] : x_j\leq x_i \}|.
$$
Equivalently, $\red(w)$ is the permutation in $\sym_n$ which is
order equivalent to $w$. For example, $\red(2659) = 1324$.

We may regard a \emph{pattern} as a function from $\sym_n$ to the set
$\N$ of natural numbers. The patterns of main interest to us are
defined as follows. Let $xyz\in\sym_3$ and $\pi = a_1 a_2 \cdots
a_n\in\sym_n$, then
$$
  (x\dd yz)\,\pi \,=\,
  |\{ a_i a_j a_{j+1} : \red(a_i a_j a_{j+1}) = xyz,
                        1\leq i < j < n \}|
$$
and similarly $(xy\dd z)\,\pi = (z\dd yx)\,\pi^r$. For instance
$$(\axbc)\, 491273865 = |\{127,138,238\}| = 3.$$
A pattern
$p=p_1\dd p_2\dd\cdots\dd p_k$ containing exactly $k-1$ dashes
is said to be of type $(|p_1|,|p_2|,\ldots,|p_k|)$.  For
example, the pattern $142\dd 5\dd 367$ is of type $(3,1,3)$, and any
classical pattern of length $k$ is of type
$(\underbrace{1,1,\ldots,1}_{k})$.

We say that a
permutation $\pi$ \emph{avoids} a pattern $p$ if $p\,\pi = 0$. The set
of all permutations in $\sym_n$ that avoids $p$ is denoted $\sym_n(p)$
and, more generally, $\sym_n(P) = \bigcap_{p\in P} \sym_n(p)$ and
$\sym(P) = \bigcup_{n \geq 0} \sym_n(P)$.

We extend the definition of reverse and complement to patterns the
following way. Let us call $\pi$ the \emph{underlying permutation} of
the pattern $p$ if $\pi$ is obtained from $p$ by deleting all the
dashes in $p$. If $p$ is a pattern with underlying permutation $\pi$,
then $p^c$ is the pattern with underlying permutation $\pi^c$ and with
dashes at precisely the same positions as there are dashes in $p$. We
define $p^r$ as the pattern we get from regarding $p$ as a word and
reading it backwards. For example, $(\axbc)^c = \cxba$ and
$(\axbc)^r = \cbxa$. Observe that
\begin{align*}
  \sigma\in\sym_n(p) &\iff \sigma^r\in\sym_n(p^r) \\
  \sigma\in\sym_n(p) &\iff \sigma^c\in\sym_n(p^c).
\end{align*}
These observations of course generalize to $\sym_n(P)$ for any set of
patterns $P$.

The operations reverse and complement generates the dihedral group
$D_2$ (the symmetry group of a rectangle).  The orbits of $D_2$ in the
set of patterns of type $(1,2)$ or $(2,1)$ will be called
\emph{symmetry classes}. For instance, the symmetry class of $\axbc$
is
$$
\{\axbc, \cxba, \abxc, \cbxa\}.
$$
We also talk about symmetry classes of sets of patterns (defined in
the obvious way). For example, the symmetry class of $\{\axbc,\cxba\}$
is
$$\{\{\axbc,\cxba\}, \{\cbxa,\abxc\}\}.$$
A set of patterns $P$ such that if $p,p'\in P$ then, for each $n$,
$|\sym_n(p)|=|\sym_n(p')|$ is called a \emph{Wilf-class}. For
instance, by Proposition~\ref{claesson1}, the Wilf-class of $\axbc$ is
$$
\{\axbc, \cxba, \abxc, \cbxa,\axcb, \cxab, \baxc, \bcxa\}.
$$
We also talk about Wilf-classes of sets of patterns (defined in the
obvious way). It is clear that symmetry classes are Wilf-classes, but
as we have seen the converse does not hold in general.

In what follows we will frequently use the well known
bijection between increasing binary trees and permutations (e.g. see
\cite[p. 24]{St97}). Let $\pi$ be any word on the alphabet
$\{1,2,3,4,\ldots\}$ with no repeated letters. If $\pi\neq\emptyword$
then we can factor $\pi$ as $\pi = \sigma\,\hat{0}\,\tau$, where
$\hat{0}$ is the minimal element of $\pi$.
Define $T(\emptyword) = \bullet$ (a leaf) and
$$T(\pi) \;=\;
\vcenter{\bintree{
            & \hat{0}\LR \\
  T(\sigma) &              & T(\tau)
}}
$$
In addition, we define $U(t)$ as the
unlabelled counterpart of the labelled tree $t$. For instance
\begin{center}
$T(316452)\;=\;
\vcenter{\bintree{
    & 1\LR \\
  3 &      & 2\LEFT \\
    & 4\LR \\
  6 &      & 5
}}
$\quad\quad\quad
$U\circ T(316452)\;=\;
\vcenter{\bintree{
          & {\Circ}\LR \\
  {\Circ} &            & {\Circ}\LEFT \\
          & {\Circ}\LR \\
  {\Circ} &            & \Circ
}}
$
\end{center}
Note that, for sake of simplicity, the leafs are not displayed.

\section{Single patterns}\label{single}

There are $3$ symmetry classes and $2$ Wilf-classes of single
patterns. The details are as follows.

\begin{proposition}[Claesson~\cite{Cl01}]
  Let $n\in\N$. We have
  $$
  |\sym(p)| =
  \begin{cases}
    B_n & \text{if }\,
    p\in\{
    \axbc, \cxba, \abxc, \cbxa
    \},\\
    B_n & \text{if }\,
    p\in\{
    \axcb, \cxab, \baxc, \bcxa
    \},\\
    C_n &\text{if }\,
    p\in\{
    \bxac, \bxca, \acxb, \caxb
    \},\\
  \end{cases}
  $$
  where $B_n$ and $C_n$ are the $n$th Bell and Catalan numbers,
  respectively.
\end{proposition}

\begin{proof}[Proof of the first case]
  Note that
  $$
  \sigma 1 \tau \in \sym(\axbc) \iff
  \begin{cases}
    \red(\sigma)\in\sym(\axbc) \\
    \red(\tau)\in\sym(\ab) \\
    \sigma 1\tau \in \sym
  \end{cases}
  $$
  where of course $\sym(\ab) =
  \{\emptyword,1,21,321,4321,\ldots\}$.  This enables us to give a
  bijection $\Phi$ between $\sym_n(\axbc)$ and the set of partitions
  of $[n]$, by induction.  Let $\Phi(\emptyword)$ be the empty
  partition. Let the first block of $\Phi(\sigma 1\tau)$ be the set of
  letters of $1\tau$, and let the rest of the blocks of $\Phi(\sigma
  1\tau)$ be as in $\Phi(\sigma)$.
\end{proof}

The most transparent way to see the above correspondence is perhaps to
view the permutation as an increasing binary tree. For instance, the tree
$$T(649752183) \;=\;
  \vcenter{\bintree{
      &        &        & 1\LR \\
      &        & 2\LEFT &         & 3\LEFT \\
      & 4\LR   &        & 8 \\
    6 &        & 5\LEFT \\
      & 7\LEFT \\
    9
    }}
$$
corresponds to the partition
$\{\{1,3,8\},\{2\},\{4,5,7,9\},\{6\}\}$.

\begin{proof}[Proof of the second case]
  This case is analogous to the previous one. We have
  $$
  \sigma 1 \tau \in \sym(\axcb) \iff
  \begin{cases}
    \red(\sigma)\in\sym(\axcb) \\
    \red(\tau)\in\sym(\ba) \\
    \sigma 1\tau \in \sym
  \end{cases}
  $$
  We give a bijection $\Phi$ between $\sym_n(\axbc)$ and the set of
  partitions of $[n]$, by induction. Let $\Phi(\emptyword)$ be the
  empty partition. Let the first block of $\Phi(\sigma 1\tau)$ be the
  set of letters of $1\tau$, and let the rest of the blocks of
  $\Phi(\sigma 1\tau)$ be as in $\Phi(\sigma)$.
\end{proof}

As an example, the tree
$$T(645792138) \;=\;
  \vcenter{\bintree{
      &      &         & 1\LR \\
      &      & 2\LEFT  &         & 3\RIGHT \\
      & 4\LR &         &         &         & 8 \\
    6 &      & 5\RIGHT \\
      &      &         & 7\RIGHT \\
      &      &         &         & 9
    }}
$$
corresponds to the partition
$\{\{1,3,8\},\{2\},\{4,5,7,9\},\{6\}\}$.

Now that we have seen the structure of $\sym(\axbc)$ and
$\sym(\axcb)$, it is trivial to give a bijection between the two
sets. Indeed, if $\Theta:\sym(\axbc)\to\sym(\axcb)$\label{Theta}
is given by $\Theta(\emptyword) = \emptyword$ and $\Theta(\sigma
1\tau) = \Theta(\sigma)\,1\,\tau^r$ then $\Theta$ is such a
bijection. Actually $\Theta$ is its own inverse.

\begin{proof}[Proof of the third case]
  It is plain that a permutation avoids $\bxac$ if and only if it
  avoids $\bxaxc$ (see \cite{Cl01}). Note that
  $$
  \sigma 1 \tau \in \sym(\bxaxc) \iff
  \begin{cases}
    \red(\sigma),\red(\tau)\in\sym(\bxaxc) \\
    \tau > \sigma\\
    \sigma 1\tau \in \sym
  \end{cases}
  $$
  where $\tau > \sigma$ means that any letter of $\tau$ is greater
  than any letter of $\sigma$. Hence we get a unique labelling of the
  binary tree corresponding to $\sigma 1 \tau$, that is, if $\pi_1,
  \pi_2\in \sym(\bxaxc)$ and $U\circ T(\pi_1) = U\circ T(\pi_2)$ then
  $\pi_1 = \pi_2$. It is well known that there are exactly $C_n$
  (unlabelled) binary trees with $n$ (internal) nodes. The validity of
  the last statement can be easily deduced from the following simple
  bijection between Dyck words and binary trees. Fixing notation, we
  let the set of Dyck words be the smallest set of words over
  $\{u,d\}$ that contains the empty word and is closed under
  $(\alpha,\beta)\mapsto u\alpha d\beta$. Now the promised bijection
  is given by $\Psi(\bullet) = \emptyword$ and
  $$
  \Psi\biggl(\vcenter{\bintree{ &\Circ\LR \\ L & & R }}\biggr)
  = u \Psi(L) d \Psi(R).
  $$
\end{proof}

\section{Pairs of patterns}\label{double}

There are $\binom {12} 2 = 66$ pairs of patterns altogether. It turns
out that there are $21$ symmetry classes and $10$ Wilf-classes. The
details are given in Table~\ref{tab}, and the numbering of the
symmetry classes in the titles of the subsections below is taken from
that table.

\begin{table}
$$
\begin{array}{|c|c|c|}
  \hline
  & \{p,q\} & |\sym_n(p,q)|\\
  \hline\hline
  1  &
  \begin{array}{c}
    \axbc, \cbxa \\
    \cxba, \abxc \\
  \end{array}
  & 0 \\
  \hline
  2  &
  \begin{array}{c}
    \axbc, \cxba \\
    \cbxa, \abxc \\
  \end{array}
  & 2(n-1) \\
  \hline
  3  &
  \begin{array}{c}
    \axbc, \bxca \\
    \cxba, \bxac \\
    \abxc, \caxb \\
    \cbxa, \acxb \\
  \end{array}
  & \binom n 2 + 1\\
  \hline
  4a &
  \begin{array}{c}
    \axbc, \bxac \\
    \cxba, \bxca \\
    \abxc, \acxb \\
    \cbxa, \caxb \\
  \end{array}
  & 2^{n-1}\\
  \hline
  4b &
  \begin{array}{c}
    \axbc, \bcxa \\
    \cxba, \baxc \\
    \abxc, \cxab \\
    \cbxa, \axcb \\
  \end{array}
  & 2^{n-1}\\
  \hline
  4c &
  \begin{array}{c}
    \axbc, \caxb \\
    \cxba, \acxb \\
    \abxc, \bxca \\
    \cbxa, \bxac \\
  \end{array}
  & 2^{n-1}\\
  \hline
  4d &
  \begin{array}{c}
    \axcb, \bxac \\
    \cxab, \bxca \\
    \acxb, \baxc \\
    \bcxa, \caxb \\
  \end{array}
  & 2^{n-1}\\
  \hline
  4e &
  \begin{array}{c}
    \axcb, \bxca \\
    \cxab, \bxac \\
    \caxb, \baxc \\
    \bcxa, \acxb \\
  \end{array}
  & 2^{n-1}\\
  \hline
  4f &
  \begin{array}{c}
    \axcb, \cxab \\
    \bcxa, \baxc \\
  \end{array}
  & 2^{n-1}\\
  \hline
  4g &
  \begin{array}{c}
    \axcb, \bcxa \\
    \cxab, \baxc \\
  \end{array}
  & 2^{n-1}\\
  \hline
\end{array}
\quad
\begin{array}{|c|c|c|}
  \hline
  & \{p,q\} & |\sym_n(p,q)|\\
  \hline\hline
  4h &
  \begin{array}{c}
    \axcb, \caxb \\
    \cxab, \acxb \\
    \baxc, \bxca \\
    \bcxa, \bxac \\
  \end{array}
  & 2^{n-1}\\
  \hline
  4i &
  \begin{array}{c}
    \bxac, \bxca \\
    \caxb, \acxb \\
  \end{array}
  & 2^{n-1}\\
  \hline
  4j &
  \begin{array}{c}
    \bxac, \acxb \\
    \bxca, \caxb \\
  \end{array}
  & 2^{n-1}\\
  \hline
  4k &
  \begin{array}{c}
    \bxac, \caxb \\
    \bxca, \acxb \\
  \end{array}
  & 2^{n-1}\\
  \hline
  5a &
  \begin{array}{c}
    \axbc, \acxb \\
    \cxba, \caxb \\
    \abxc, \bxac \\
    \cbxa, \bxca \\
  \end{array}
  &
  \begin{array}{c}
  M_n \\
  (\text{Motzkin no.})\\
  \end{array}\\
  \hline
  5b &
  \begin{array}{c}
    \axbc, \baxc \\
    \cxba, \bcxa \\
    \abxc, \axcb \\
    \cbxa, \cxab \\
  \end{array}
  &
  \begin{array}{c}
  M_n \\
  (\text{Motzkin no.})\\
  \end{array}\\
  \hline
  6 &
  \begin{array}{c}
    \axcb, \baxc \\
    \cxab, \bcxa \\
  \end{array}
  & a_n \\
  \hline
  7 &
  \begin{array}{c}
    \axbc, \cxab \\
    \cxba, \axcb \\
    \bcxa, \abxc \\
    \cbxa, \baxc \\
  \end{array}
  & b_n \\
  \hline
  8 &
  \begin{array}{c}
    \axbc, \axcb \\
    \cxba, \cxab \\
    \baxc, \abxc \\
    \cbxa, \bcxa \\
  \end{array}
  &
  \begin{array}{c}
    I_n \\
    (\text{\# involutions})\\
  \end{array}\\
  \hline
  9 &
  \begin{array}{c}
    \axcb, \acxb \\
    \cxab, \caxb \\
    \baxc, \bxac \\
    \bcxa, \bxca \\
  \end{array}
  &
  \begin{array}{c}
    C_n \\
   (\text{Catalan no.})\\
  \end{array}\\
  \hline
  10 &
  \begin{array}{c}
    \axbc, \abxc \\
    \cxba, \cbxa \\
  \end{array}
  &  B_n^* \;\;(\text{Bessel no.})\\
  \hline
\end{array}
$$
\caption{}\label{tab}
\end{table}

\subsection*{Symmetry class 1}

We have
$$
\sigma 1 \tau \in \sym(\axbc,\cbxa) \iff
\begin{cases}
  \red(\sigma)\in\sym(\ba,\axbc) \\
  \red(\tau)\in\sym(\ab,\cbxa) \\
  \sigma 1\tau \in \sym
\end{cases}
$$
The result now follows from $\sym(\ba,\axbc)= \{\emptyword,1,12\}$ and
$\sym(\ab,\cbxa) = \{\emptyword,1,21\}$.

\subsection*{Symmetry class 2}

Since $\cxba$ is the complement of $\axbc$, the cardinality of
$\sym_n(\axbc,\cxba)$ is twice the number of
permutations in $\sym_n(\axbc,\cxba)$ in which $1$ precedes $n$. In
addition, $1$ and $n$ must be adjacent letters in a permutation
avoiding $\axbc$ and $\cxba$. Let $\sigma 1 n \tau$ be such a
permutation. Note that $\tau$ must be both increasing and
decreasing, that is, $\tau\in\{\emptyword,2,3,4,\ldots,n-1\}$, so
there are $n-1$ choices for $\tau$. Furthermore, there is exactly one
permutation in $\sym_n(\axbc,\cxba)$ of the form $\sigma 1 n$,
namely $(\lceil\frac {n+1} 2\rceil,\ldots,n-2,3,n-1,2,n,1)$, and
similarly there is exactly one of the form $\sigma 1 n k$ for each
$k\in\{2,3,\ldots,n-1\}$. This completes our argument.

\subsection*{Symmetry class 3}

Note that
$$
\sigma 1 \tau \in \sym(\axbc,\bxca) \iff
\begin{cases}
  \red(\sigma),\red(\tau) \in \sym(\ab)\\
  \sigma 1\tau \in \sym(\bxca)
\end{cases}
$$
It is now rather easy to see that $\pi\in\sym_n(\axbc,\bxca)$ if
and only if $\pi = n \cdots 2 1$ or $\pi$ is constructed in the
following way. Choose $i$ and $j$ such that $1\leq j < i \leq n$.
Let $\pi(i-1) = 1$, $\pi(i) = n+1-j$ and arrange the rest of the
elements so that $\pi(1)>\pi(2)>\cdots >\pi(i-1)$ and
$\pi(i)>\pi(i+1)>\cdots>\pi(n)$ (this arrangement is unique). Since
there are $\binom n 2$ ways of choosing $i$ and $j$ we get the
desired result.

\subsection*{Symmetry class 4a}

We have
$$
\sigma 1 \tau \in \sym(\axbc,\bxac) \iff
\begin{cases}
  \red(\sigma)\in\sym(\axbc,\bxac) \\
  \red(\tau)\in\sym(\ab) \\
  \sigma > \tau \\
  \sigma 1\tau \in \sym,
\end{cases}
$$
where $\sigma > \tau$ means that any letter of $\tau$ is greater
than any letter of $\sigma$. This enables us to give a bijection
between $\sym_n(\axbc,\bxac)$ and the set of compositions (ordered
formal sums) of $n$. Indeed, such a bijection $\Psi$ is given by
$\Psi(\emptyword) = \emptyword$ and $\Psi( \sigma 1 \tau) =
\Psi(\sigma) + |1\tau|$.

As an example, the tree
$$U\circ T(958764132) \;=\;
  \vcenter{\ulbintree{
          &            &            & \circ\LR \\
          &            & \circ\LEFT &           & \circ\LEFT \\
          & \circ\LR   &            & \circ \\
    \circ &            & \circ\LEFT \\
          & \circ\LEFT \\
    \circ
    }}
$$
corresponds to the composition $1 + 4 + 1 + 3$ of $9$.

\subsection*{Symmetry class 4b}

We have
$$
  \sigma 1 \tau \in \sym(\axbc,\bcxa) \iff
  \begin{cases}
    \red(\sigma),\red(\tau)\in\sym(\ab)\\
    \sigma 1\tau \in \sym
  \end{cases}
$$
Hence a permutation in $\sym(\axbc,\bcxa)$ is given by the
following procedure. Choose a subset $S\subseteq
\{2,3,4,\ldots,n\}$, let $\sigma$ be the word obtained by writing
the elements of $S$ in decreasing order, and let $\tau$ be the word
obtained by writing the elements of $\{2,3,4,\ldots,n\}\setminus S$
in decreasing order.

For instance, the tree
$$ T(421653) \;=\; \vcenter{\bintree{
       &        & 1\LR \\
       & 2\LEFT &        & 3\LEFT \\
     4 &        & 5\LEFT \\
       & 6 \\
    }}
$$
corresponds to the subset $\{2,4\}$ of $\{2,3,4,5,6\}$.

\subsection*{Symmetry class 4c}

This case is essentially identical to the case dealt with in (4a).

\subsection*{Symmetry class 4c}

The bijection $\Theta$ between $\sym(\axbc)$ and $\sym(\axcb)$ (see
page \pageref{Theta}) provides a one-to-one correspondence between
$\sym_n(\axcb,\bxac)$ and $\sym_n(\axbc,\bxac)$. Consequently the
result follows from (4a).

\subsection*{Symmetry class 4e}

We have
$$
\sigma 1 \tau \in \sym(\cxab,\bxac) \iff
\begin{cases}
  \red(\sigma),\red(\tau)\in\sym(\cxab,\bxac) \\
  \sigma = \emptyword \;\text{ or }\; \tau = \emptyword \\
  \sigma 1\tau \in \sym
\end{cases}
$$
Thus a bijection between $\sym_n(\cxab,\bxac)$ and
$\{0,1\}^{n-1}$ is given by $\Psi(\emptyword) = \emptyword$ and
$$
\Psi(\sigma 1 \tau) = x\Psi(\sigma\tau)
\text{ where } x =
\begin{cases}
  1 & \text{if } \sigma \neq \emptyword,\\
  0 & \text{if } \tau \neq \emptyword,\\
  \emptyword & \text{ otherwise.}
\end{cases}
$$

As an example, the tree
$$U\circ T(136542) \;=\; \vcenter{\ulbintree{
    & \circ\RIGHT \\
    &             & \circ\LEFT \\
    & \circ\RIGHT \\
    &             & \circ\LEFT \\
    & \circ\LEFT \\
    \circ
    }}
$$
corresponds to $01011\in\{0,1\}^5$.

\subsection*{Symmetry class 4f}

Since $\cxab$ is the complement of $\axcb$, the cardinality
of $\sym_n(\axcb,\cxab)$ is twice the number of
permutations in $\sym_n(\axcb,\cxab)$ in which $1$ precedes $n$.  In
addition, $n$ must be the last letter in such a permutation or else
a hit of $\axcb$ would be formed. We have
\begin{eqnarray*}
  \sigma 1 \tau n \in \sym(\axcb,\cxab) & \iff &
  \begin{cases}
    \red(\sigma 1 \tau)\in\sym(\axcb,\cxab) \\
    \red(\tau) \in \sym(\ba) \\
    \sigma 1\tau \in \sym
  \end{cases} \\
  &\iff&
  \begin{cases}
    \red(\sigma)\in\sym(\axcb,\cxab) \\
    \red(\tau) \in \sym(\ba) \\
    \sigma < \tau \\
    \sigma 1\tau \in \sym
  \end{cases} \\
\end{eqnarray*}
The rest of the proof follows the same lines as the proof of (4a).

\subsection*{Symmetry class 4g}

We can copy almost verbatim the proof of (4e);
indeed, it is easy to see that
$\sym_n(\axcb,\bcxa)=\sym_n(\axcb,\bxca)$.

\subsection*{Symmetry class 4h}

We can copy almost verbatim the proof of (4f);
indeed, it is easy to see that
$\sym_n(\axcb,\caxb)=\sym_n(\axcb,\cxab)$.

\subsection*{Symmetry class 4i}

$|\sym_n(\bxac,\bxca)|=|\sym_n(\bxaxc,\bxcxa)| = 2^{n-1}$ by
\cite[Lemma $5$(d)]{SiSc85}.

\subsection*{Symmetry class 4j}

$|\sym_n(\bxac,\acxb)|=|\sym_n(\axcxb,\bxaxc)| =
2^{n-1}$ by \cite[Lemma $5$(b)]{SiSc85}.

\subsection*{Symmetry class 4k}

$|\sym_n(\bxac,\caxb)|=|\sym_n(\bxaxc,\cxaxb)| =
2^{n-1}$ by \cite[Lemma $5$(c)]{SiSc85}.

\subsection*{Symmetry class 5a}

See Proposition~\ref{claesson2}.

\subsection*{Symmetry class 5b}

We give a bijection $$\Lambda:\sym_n(\axbc, \baxc)\to\sym_n(\axbc,
\acxb)$$ by means of induction. Let $\pi\in\sym_n(\axbc, \baxc)$.
Define $\Lambda(\pi) = \pi$ for $n\leq 1$. Assume $n\geq 2$ and
$\pi=a_1 a_2 \cdots a_n$. It is plain that either $a_1 = n$ or $a_2
= n$, so we can define $\Lambda(\pi)$ by
$$
\begin{cases}
  (a'_{1}+1,\ldots,a'_{n-1}+1,a'_{n-2}+1,1)
  & \text{if }
  \begin{cases}
    a_1 = n \;\;\;\text{ and }\\
    a'_1\cdots a'_{n-1} =\Lambda(a_{2}a_{3}a_{4}\cdots a_{n}),
  \end{cases}\\
  (a'_{1}+1,\ldots,a'_{n-1}+1,1,a'_{n-2}+1)
  & \text{if }
  \begin{cases}
    a_2 = n \;\;\;\text{ and }\\
    a'_1\cdots a'_{n-1} = \Lambda(a_{1}a_{3}a_{4}\cdots a_{n}).
  \end{cases}
\end{cases}
$$
Observing that if $\sigma \in \sym_n(\axbc,\acxb)$ then
$\sigma(n-1)=1$ or $\sigma(n)=1$, it easy to find the inverse of
$\Lambda$.

\subsection*{Symmetry class 6}

In \cite{Cl01} Claesson introduced the notion of a monotone
partition. A partition is \emph{monotone} if its non-singleton
blocks can be written in increasing order of their least element
and increasing order of their greatest element, simultaneously. He
then proved that monotone partitions and non-overlapping
partitions are in one-to-one correspondence. Non-overlapping
partitions were first studied by Flajolet and Schot in
\cite{FlSc90}. A partition $\pi$ is \emph{non-overlapping} if for
no two blocks $A$ and $B$ of $\pi$ we have $\min A < \min B < \max
A < \max B$. Let $B^*_n$ be the number of non-overlapping
partitions of $[n]$; this number is called the $n$th \emph{Bessel
number}. Proposition~\ref{claesson2} tells us that there is a
bijection between non-overlapping partitions and permutations
avoiding $\axbc$ and $\abxc$. Below we define a new class of
partitions called strongly monotone partitions and then show that
there is a bijection between strongly monotone partitions and
permutations avoiding $\axcb$ and $\baxc$.

\begin{definition}\label{smon}
  Let $\pi$ be an arbitrary partition whose blocks
  $\{A_1,\ldots,A_k\}$ are ordered so that for all $i\in[k-1]$, $\min
  A_i>\min A_{i+1}$. If $\max A_i>\max A_{i+1}$ for all $i\in [k-1]$,
  then we call $\pi$ a \emph{strongly monotone partition}.
\end{definition}

In other words a partition is strongly monotone if its blocks can be
written in increasing order of their least element and increasing
order of their greatest element, simultaneously. Let us denote by
$a_n$ the number of strongly monotone partitions of $[n]$.
The sequence
$\{a_n\}_0^{\infty}$ starts with
$$
1, 1, 2, 4, 9, 22, 58, 164, 496, 1601, 5502, 20075, 77531, 315947, 1354279.
$$
It is routine to derive the continued fraction expansion
$$
\sum_{n\geq 0} a_n x^n =
\cfrac{1}{1 - 1 \cdot x -
  \cfrac{x^2}{1 -  1 \cdot x -
    \cfrac{x^2}{1 -  2 \cdot x -
      \cfrac{x^2}{1 -  3 \cdot x -
        \cfrac{x^2}{1 -  4 \cdot x -
          \cfrac{x^2}{\quad\ddots
            }}}}}}
$$
using the standard machinery of Flajolet~\cite{Fl80} and
Fran{\c{c}}on and Viennot~\cite{FrVi79}. One can also note that there
is a one-to-one correspondence between strongly monotone partitions
and non-overlapping partition, $\pi$, such that if $\{x\}$ and $B$ are
blocks of $\pi$ then either $x<\min B$ or $\max B < x$. In addition,
we observe that
$$ \sum_{n\geq 0} a_n x^n = \frac 1 {1-x-x^2B^*(x)},$$
where $B^*(x) = \sum_{n\geq 0} B^*_n x^n$ is the ordinary generating
function for the Bessel numbers.

Suppose $\pi\in\sym_n$ has $k+1$ left-to-right minima
$1,1',1'',\ldots,1^{(k)}$ such that
$$
1<1'<1''<\cdots<1^{(k)} \text{, and }
\pi = 1^{(k)}\tau^{(k)} \cdots 1'\tau' 1 \tau.
$$
Then $\pi$ avoids $\axcb$ if and only if, for each
$i$, $\tau^{(i)}\in\sym(\ba)$. If $\pi$ avoids $\axcb$ and $x_i =
\max 1^{(i)}\tau^{(i)}$ then $\pi$ avoids $\baxc$ precisely when
$x_0 < x_1 < \cdots < x_k$. This follows from observing that the
only potential $(\baxc)$-subwords of $\pi$ are $x_{i+1}1^{(k)}x_j$
with $j\leq i$.

Mapping $\pi$ to the partition $ \{ 1\sigma, 1'\sigma', \ldots,
1^{(k)}\tau^{(k)} \}$ we thus get a one-to-one correspondence
between permutations in $\sym_n(\axcb,\baxc)$ and strongly monotone
partitions of $[n]$.

\subsection*{Symmetry class 7}

Let the sequence $\{b_n\}$ be defined by $b_0=1$ and, for $n \geq -2$,
$$b_{n+2} = b_{n+1} + \sum_{k=0}^n\binom n k b_k.
$$
The first few of the numbers $b_n$ are
$$
1,1,2,4,9,23,65,199,654,2296,\ldots
$$

Suppose $\pi\in\sym_n$ has $k+1$ left-to-right minima
$1,1',1'',\ldots,1^{(k)}$ such that
$$
1<1'<1''<\cdots<1^{(k)} \text{, and }
\pi = 1^{(k)}\tau^{(k)} \cdots 1'\tau' 1 \tau.
$$
Then $\pi$ avoids $\axbc$ if and only if, for each
$i$, $\tau^{(i)}\in\sym(\ab)$. If $\pi$ avoids $\axbc$ and $x_i =
\max 1^{(i)}\tau^{(i)}$ then $\pi$ avoids $\cxab$ precisely when
$$ j>i \;\text{ and }\, x_i \neq 1^{(i)}\;\;\Longrightarrow\;\;x_j < x_i. $$
This follows
from observing that the only potential $(\cxab)$-subwords of $\pi$
are $x_j1^{(k)}x_i$ with $j\leq i$. Thus we have established
$$
\sigma 1 \tau \in \sym_n(\axbc, \cxab) \iff
\begin{cases}
  \red(\sigma)\in\sym(\axbc, \cxab)\\
  \tau\neq\emptyword \;\;\Rightarrow\;\;
  \tau = \tau' n \,\text{ and } \red(\tau')\in\sym(\ab)\\
  \sigma 1\tau \in \sym_n
\end{cases}
$$
If we know that $\sigma 1 \tau' n \in \sym_{n}(\axbc, \cxab)$
and $\red(\tau')\in\sym_k(\ab)$ then there are $\binom {n-2} k$
candidates for $\tau'$. In this way the recursion follows.

\subsection*{Symmetry class 8}

See Proposition~\ref{claesson2}.

\subsection*{Symmetry class 9}

$\sym_n(\axcb,\acxb)=\sym_n(\axcxb)$.

\subsection*{Symmetry class 10}

See Proposition~\ref{claesson2}.

\section{More than two patterns}\label{multi}
\newcommand{\X}{\times}

Let $P$ be a set of patterns of type $(1,2)$ or $(2,1)$. With the
aid of a computer we have calculated the cardinality of
$\sym_n(P)$ for sets $P$ of three or more patterns. From these
results we arrived at the plausible conjectures of table~\ref{tb1}
(some of which are trivially true).  We use the notation
\mbox{$m\X n$} to express that there are $m$ symmetric classes
each of which contains $n$ sets. Moreover, we denote by $F_n$ the
$n$th \emph{Fibonacci number} ($F_0=F_1=1, F_{n+1}=F_n+F_{n-1}$).

\begin{table}
  \begin{tabular}{|c|c|}
    \hline
    \begin{minipage}{38ex}
      \vspace{1ex}
      For $|P|=3$ there are $220$ sets, $55$ \mbox{symmetry}
      classes and $9$ Wilf-classes.
      $$
      \begin{array}{rr}
        \text{cardinality}   & \# \text{ sets} \\
        \hline
        0                    & 7\X 4  \\
        3                    & 1\X 4  \\
        n                    & 24\X 4 \\
        1+\binom{n}{2}       & 2\X 4  \\
        F_n                  & 7\X 4  \\
        \binom{n}{[n/2]}     & 1\X 4  \\
        2^{n-2}+1            & 1\X 4  \\
        2^{n-1}              & 10\X 4 \\
        M_n                  & 2\X 4  \\
      \end{array}
      $$
      \vspace{.3ex}
    \end{minipage}
    &
    \begin{minipage}{38ex}
      \vspace{1ex}
      For $|P|=4$ there are $495$ sets, $135$ \mbox{symmetry}
      classes, and $9$ Wilf-classes.
      $$
      \begin{array}{rr}
        \text{cardinality}   & \# \text{ sets} \\
        \hline
        0                    & 1\X 1+6\X 2+30\X 4  \\
        2                    & 2\X 1+5\X 2+35\X 4  \\
        3                    & 1\X 4               \\
        n                    & 37\X 4+1\X 2        \\
        1+\binom{n}{2}       & 1\X 4               \\
        F_n                  & 9\X 4+1\X 2         \\
        \binom{n}{[n/2]}     & 1\X 2               \\
        2^{n-2}+1            & 1\X 2               \\
        2^{n-1}              & 1\X 4+3\X 2         \\
      \end{array}
      $$
      \vspace{.3ex}
    \end{minipage}
    \\ \hline
    \begin{minipage}{38ex}
      \vspace{1ex}
      For $|P|=5$ there are $792$ sets, $198$ \mbox{symmetry}
      classes, and $5$ Wilf-classes.
      $$
      \begin{array}{rr}
        \text{cardinality}   & \# \text{ sets} \\
        \hline
        0                    & 84\X 4 \\
        1                    & 16\X 4 \\
        2                    & 74\X 4 \\
        n                    & 20\X 4 \\
        F_n                  & 4\X 4  \\
      \end{array}
      $$
      \vspace{.3ex}
    \end{minipage}
    &
    \begin{minipage}{38ex}
      \vspace{1ex}
      For $|P|=6$ there are $924$ sets, $246$ \mbox{symmetry}
      classes, and $4$ Wilf-classes.
      $$
      \begin{array}{rr}
        \text{cardinality}   & \# \text{ sets} \\
        \hline
        0                    & 17\X 2+124\X 4 \\
        1                    & 4\X 2+38\X 4   \\
        2                    & 7\X 2+51\X 4   \\
        n                    & 1\X 2+3\X 4    \\
        F_n                  & 1\X 2          \\
      \end{array}
      $$
      \vspace{.3ex}
    \end{minipage}
    \\ \hline
    \begin{minipage}{38ex}
      \vspace{1ex}
      For $|P|=7$ there are $792$ sets, $198$ \mbox{symmetry}
      classes, and $3$ Wilf-classes.
      $$
      \begin{array}{rr}
        \text{cardinality}   & \# \text{ sets} \\
        \hline
        0                    & 140\X 4\\
        1                    & 40\X 4 \\
        2                    & 18\X 4 \\
      \end{array}
      $$
      \vspace{.3ex}
    \end{minipage}
    &
    \begin{minipage}{38ex}
      \vspace{1ex}
      For $|P|=8$ there are $495$ sets, $135$ \mbox{symmetry}
      classes, and $3$ Wilf-classes.
      $$
      \begin{array}{rr}
        \text{cardinality}   & \# \text{ sets} \\
        \hline
        0                    & 2\X 1+14\X 2+94\X 4 \\
        1                    & 4\X 2+18\X 4        \\
        2                    & 1\X 1+2\X 4         \\
      \end{array}
      $$
      \vspace{.3ex}
    \end{minipage}
    \\ \hline
    \begin{minipage}{38ex}
      \vspace{1ex}
      For $|P|=9$ there are $220$ sets, $55$ \mbox{symmetry}
      classes, and $2$ Wilf-classes.
      $$
      \begin{array}{rr}
        \text{cardinality}   & \# \text{ sets} \\
        \hline
        0                    &  50\X 4      \\
        1                    &  5\X 4       \\
      \end{array}
      $$
      \vspace{.3ex}
    \end{minipage}
    &
    \begin{minipage}{38ex}
      \vspace{1ex}
      For $|P|=10$ there are $66$ sets, $21$ \mbox{symmetry}
      classes, and $2$ Wilf-classes.
      $$
      \begin{array}{rr}
        \text{cardinality}   & \# \text{ sets} \\
        \hline
        0                    &  8\X 2+12\X 4  \\
        1                    &  1\X 2         \\
      \end{array}
      $$
      \vspace{.3ex}
    \end{minipage}
    \\ \hline
    \begin{minipage}{38ex}
      \vspace{1ex}
      For $|P|=11$ there are $12$ sets, $3$ \mbox{symmetry}
      classes, and $1$ Wilf-class.
      $$
      \begin{array}{rr}
        \text{cardinality}   & \# \text{ sets} \\
        \hline
        0                    & 3\X 4 \\

      \end{array}
      $$
      \vspace{.3ex}
    \end{minipage}
    &
    \begin{minipage}{38ex}
      \vspace{1ex}
      For $|P|=12$ there is $1$ set, $1$ \mbox{symmetry}
      class, and $1$ Wilf-class.
      $$
      \begin{array}{rr}
        \text{cardinality}   & \# \text{ sets} \\
        \hline
        0                    & 1\X 1 \\
      \end{array}
      $$
      \vspace{.3ex}
    \end{minipage}
    \\
    \hline
  \end{tabular}
  \caption{The cardinality of $\sym_n(P)$ for $|P|>2$.}\label{tb1}
\end{table}

\section*{Acknowledgements}
The first author wishes to express his gratitude towards Einar
Steingr\'{\i}- msson, Kimmo Eriksson, and Mireille Bousquet-M\'{e}lou;
Einar for his guidance and infectious enthusiasm; Kimmo for useful
suggestions and a very constructive discussion on the results of this
paper; Mireille for her great hospitality during a stay at LaBRI,
where some of the work on this paper was done.\\
We would like to thank N. J. A. Sloane for his excellent web site
``The On-Line Encyclopedia of Integer Sequences''
\begin{center}
{\tt http://www.research.att.com/\~{}njas/sequences/}.
\end{center}
It is simply an indispensable tool for all studies concerned with
integer sequences.

\bibliographystyle{plain}

\end{document}